\input amstex
\input amsppt.sty
\magnification=\magstep1
\hsize=30truecc
\vsize=22.2truecm
\baselineskip=16truept
\NoBlackBoxes
\TagsOnRight \pageno=1 \nologo
\def\Z{\Bbb Z}
\def\N{\Bbb N}

\def\l{\left}
\def\r{\right}
\def\bg{\bigg}
\def\({\bg(}
\def\[{\bg\lfloor}
\def\){\bg)}
\def\]{\bg\rfloor}
\def\t{\text}
\def\f{\frac}

\def\p{\ (\roman{mod}\ p)}

\def\sm{\setminus}

\def\bi{\binom}
\def\eq{\equiv}

\def\ls{\leqslant}
\def\gs{\geqslant}
\def\mo{\roman{mod}}

\def\ve{\varepsilon}
\def\al{\alpha}

\def\Proof{\noindent{\it Proof}}

\def\Remark{\medskip\noindent{\it  Remark}}

\hbox {Sci. China Math. 57(2014), no.\,7, 1375--1400}
\medskip
\topmatter
\title Congruences involving generalized central trinomial coefficients\endtitle
\author Zhi-Wei Sun\endauthor
\leftheadtext{Zhi-Wei Sun}
\rightheadtext{Congruences involving central trinomial coefficients}
\affil Department of Mathematics, Nanjing University\\
 Nanjing 210093, People's Republic of China
  \\ E-mail: {\tt zwsun\@nju.edu.cn}
\endaffil
\abstract For integers $b$ and $c$ the generalized
central trinomial coefficient $T_n(b,c)$ denotes the coefficient of $x^n$ in the
expansion of $(x^2+bx+c)^n$.
Those $T_n=T_n(1,1)\ (n=0,1,2,\ldots)$ are the usual central trinomial coefficients,
and $T_n(3,2)$ coincides with the Delannoy number $D_n=\sum_{k=0}^n\bi nk\bi{n+k}k$ in combinatorics.
We investigate congruences involving generalized central trinomial coefficients systematically.
Here are some typical results:  For each $n=1,2,3,\ldots$
we have
$$\sum_{k=0}^{n-1}(2k+1)T_k(b,c)^2(b^2-4c)^{n-1-k}\eq0\pmod{n^2}$$  and in particular $n^2\mid\sum_{k=0}^{n-1}(2k+1)D_k^2$;
if $p$ is an odd prime then
$$\sum_{k=0}^{p-1}T_k^2\eq\l(\f{-1}p\r)\ (\mo\ p)\ \ \t{and}\ \ \sum_{k=0}^{p-1}D_k^2\eq\l(\f 2p\r)\ (\mo\ p),$$
where $(-)$ denotes the Legendre symbol.
We also raise several conjectures some of which involve parameters in the representations of primes by certain
binary quadratic forms.
\endabstract
\thanks 2010 {\it Mathematics Subject Classification}.\,Primary 11A07, 11B75;
Secondary  05A10, 05A15, 11B65, 11E25.
\newline\indent {\it Keywords}. Congruences, central trinomial coefficients, Motzkin numbers,
central Delannoy numbers.
\newline\indent Supported by the National Natural Science
Foundation (grant 11171140) of China and the PAPD of Jiangsu Higher
Education Institutions.
\endthanks
\endtopmatter
\document

\heading{1. Introduction}\endheading

For $n\in\N=\{0,1,2,\ldots\}$, the $n$th central trinomial
coefficient $$T_n=[x^n](1+x+x^2)^n$$ is the
coefficient of $x^n$ in the expansion of $(1+x+x^2)^n$. Since $T_n$ is the constant term of $(1+x+x^{-1})^n$, by the multi-nomial theorem we see that
$$T_n=\sum_{k=0}^{\lfloor n/2\rfloor}\f{n!}{k!k!(n-2k)!}=\sum_{k=0}^{\lfloor n/2\rfloor}\bi n{2k}\bi{2k}k=\sum_{k=0}^n\bi nk\bi{n-k}k.$$
Central trinomial coefficients arise naturally in enumerative combinatorics (cf. Sloane [Sl]), e.g.,
$T_n$ is the number of lattice paths from the point $(0, 0)$ to $(n, 0)$ with only allowed steps $(1,0)$, $(1, 1)$ and $(1, -1)$.
As G. E. Andrews [A] pointed out, central trinomial coefficients were first studied by L. Euler. In 1987, Andrews and R. J. Baxter [AB]
found that the $q$-analogues of central trinomial coefficients have applications in the hard hexagon model.

For $n\in\N$ the $n$th Motzkin number is defined by
$$M_n=\sum_{k=0}^{\lfloor n/2\rfloor}\bi n{2k}C_k,$$
where $C_k$ denotes the $k$th Catalan number $\f1{k+1}\bi{2k}k=\bi{2k}k-\bi{2k}{k+1}$.
It is known that $M_n$ equals the number of paths from $(0,0)$ to
$(n,0)$ which never dip below the line $y=0$ and are made up of the only allowed steps $(1,0)$, $(1,1)$ and $(1,-1)$ (cf. [Sl]).

Surprisingly we find that central trinomial coefficients and Motzkin numbers have
nice congruence properties despite their combinatorial backgrounds. For example, we have the following conjecture.
(As usual, for an integer $a$ and an odd prime $p$, the notation  $(\f ap)$ stands for the Legendre symbol.)

\proclaim{Conjecture 1.1} {\rm (i)} For any $n\in\Z^+=\{1,2,3,\ldots\}$ we have
$$\sum_{k=0}^{n-1}(8k+5)T_k^2\eq0\ (\mo\ n).$$
If $p$ is a prime, then
$$\sum_{k=0}^{p-1}(8k+5)T_k^2\eq3p\l(\f p3\r)\ (\mo\ p^2).$$

{\rm (ii)} Let $p>3$ be a prime. Then
$$\align\sum_{k=0}^{p-1}M_k^2&\eq(2-6p)\l(\f p3\r)\ (\mo\ p^2),
\\\sum_{k=0}^{p-1}kM_k^2&\eq(9p-1)\l(\f p3\r)\ (\mo\ p^2),
\\\sum_{k=0}^{p-1}M_kT_k&\eq\f 43\l(\f p3\r)+\f p6\l(1-9\l(\f p3\r)\r)\ (\mo\ p^2),
\\\sum_{k=0}^{p-1}\f{M_kT_k}{(-3)^k}&\eq\f p2\l(\l(\f p3\r)-1\r)\ (\mo\ p^2),
\endalign$$
and
$$\sum_{k=0}^{p-1}\f{T_kH_k}{3^k}\eq \f{3+(\f p3)}2-p\l(1+\l(\f p3\r)\r)\ (\mo\ p^2),$$
where $H_k$ denotes the harmonic number $\sum_{0<j\ls k}1/j$.
\endproclaim

Given $b,c\in\Z$,  we define the {\it generalized central trinomial coefficients}
$$\align T_n(b,c):=&[x^n](x^2+bx+c)^n=[x^0](b+x+cx^{-1})^n
\\=&\sum_{k=0}^{\lfloor n/2\rfloor}\bi n{2k}\bi{2k}kb^{n-2k}c^k=\sum_{k=0}^{\lfloor n/2\rfloor}\bi{n-k}k\bi nkb^{n-2k}c^k
\endalign$$
and introduce the {\it generalized Motzkin numbers}
$$M_n(b,c):=\sum_{k=0}^{\lfloor n/2\rfloor}\bi n{2k}C_kb^{n-2k}c^k=\sum_{k=0}^{\lfloor n/2\rfloor}\bi{n-k}k\bi nk\f{b^{n-2k}c^k}{k+1}$$
$(n=0,1,2,\ldots)$. Note that
$$T_n=T_n(1,1),\ M_n=M_n(1,1),$$
$$T_n(2,1)=[x^n](x+1)^{2n}=\bi{2n}n,$$
and
$$M_n(2,1)=\sum_{k=0}^{\lfloor n/2\rfloor}\bi n{2k}C_k2^{n-2k}=C_{n+1}.$$
Thus $T_n(b,c)$ can be viewed a natural common extension of central binomial coefficients and
 central trinomial coefficients, while $M_n(b,c)$ can be viewed as a natural common extension of Catalan numbers
 and Motzkin numbers.
Let $d=b^2-4c$. H. S. Wilf [W, p.\,159] observed that if $\ve>0$ is sufficiently small then
$$\sum_{n=0}^\infty T_n(b,c)x^n=\f1{\sqrt{1-2bx+dx^2}}$$
for $|x|<\ve$. This implies the recurrence
$$(n+1)T_{n+1}(b,c)=(2n+1)bT_n(b,c)-dnT_{n-1}(b,c)\ \ (n\in\Z^+).$$
(See also T. D. Noe [N].) Also, the Zeilberger algorithm (cf. [PWZ, pp.\,101--119]) yields the recursion
$$(n+3)M_{n+1}(b,c)=b(2n+3)M_n(b,c)-dnM_{n-1}(b,c)\ (n=1,2,3,\ldots)$$
which implies that
$$2cx^2\sum_{n=0}^\infty M_n(b,c)x^n=1-bx-\sqrt{1-2bx+dx^2}.$$

The central Delannoy numbers (see [CHV]) are defined by
$$D_n=\sum_{k=0}^n\bi nk\bi{n+k}k=\sum_{k=0}^n\bi {n+k}{2k}\bi{2k}k\ (n\in\N).$$
Such numbers also arise in many enumeration problems in combinatorics (cf. [Sl]); for example, $D_n$
is the number of lattice paths from the point $(0,0)$ to $(n,n)$ with steps $(1,0),(0,1)$ and $(1,1)$.
For $n\in\N$ we define the polynomial
$$D_n(x)=\sum_{k=0}^n\bi nk\bi{n+k}kx^k.$$
Note that $D_n((x-1)/2)$ coincides with the well-known Legendre polynomial $P_n(x)$ of degree $n$.
It is known that
$$\sum_{n=0}^\infty P_n(t)x^n=\f1{\sqrt{1-2tx+x^2}}.$$
Thus, if $b,c\in\Z$ and $d=b^2-4c\not=0$ then
$$\sum_{n=0}^\infty T_n(b,c)\l(\f x{\sqrt d}\r)^n=\f1{\sqrt{1-2bx/\sqrt d+d(x/\sqrt d)^2}}=\sum_{n=0}^\infty P_n\l(\f b{\sqrt d}\r)x^n$$
and hence
$$T_n(b,c)=(\sqrt d)^nP_n\l(\f b{\sqrt d}\r).$$
It follows that $T_n(2x+1,x^2+x)=P_n(2x+1)=D_n(x)$
for all $x\in\Z$; in particular, $D_n=T_n(3,2)$.

Motivated by Conjecture 1.1 we investigate congruences involving generalized central trinomial coefficients as well as generalized Motzkin
numbers.

Now we state the main results of this paper.

\proclaim{Theorem 1.2} Let $p$ be an odd prime and let $b,c\in\Z$.

{\rm (i)} For any integer $m\not\eq0\pmod p$, we have
$$\sum_{k=0}^{p-1}\f{T_k(b,c)}{m^k}\eq\l(\f{(m-b)^2-4c}p\r)\ (\mo\ p)\tag1.1$$
and
$$2c\sum_{k=0}^{p-1}\f{M_k(b,c)}{m^k}\eq(m-b)^2-((m-b)^2-4c)\l(\f{(m-b)^2-4c}p\r)\ (\mo\ p).\tag1.2$$

{\rm (ii)} If $p$ does not divide $d=b^2-4c$, then
we have
$$\sum_{k=0}^{p-1}\f{T_k(b,c)^2}{d^k}\eq\l(\f {cd}p\r)\ (\mo\ p).\tag1.3$$
If $b\not\eq 2c\ (\mo\ p)$, then
$$\sum_{k=0}^{p-1}\f{T_k(b,c^2)^2}{(b-2c)^{2k}}\eq\l(\f{-c^2}p\r)\ (\mo\ p).\tag1.4$$

{\rm (iii)} Assume that $p\nmid c$. If $d=b^2-4c\not\eq0\ (\mo\ p)$,
then
$$\sum_{k=0}^{p-1}\f{T_k(b,c)M_k(b,c)}{d^k}\eq0\ (\mo\ p).\tag1.5$$
If $D=b^2-4c^2\not\eq0\ (\mo\ p)$, then
$$\sum_{k=0}^{p-1}\f{T_k(b,c^2)M_k(b,c^2)}{(b-2c)^{2k}}\eq\f{4b}{b+2c}\l(\f Dp\r)\ (\mo\ p).\tag1.6$$
\endproclaim

{\it Example} 1.3. Let $p>3$ be a prime. Applying Theorem
1.2(ii)-(iii) with $b=c=1$ we get
$$\gather\sum_{k=0}^{p-1}\f{T_k^2}{(-3)^k}\eq\l(\f p3\r)\ (\mo\ p),\ \sum_{k=0}^{p-1}\f{T_kM_k}{(-3)^k}\eq0\pmod{p},\tag1.7
\\\sum_{k=0}^{p-1}T_k^2\eq\l(\f{-1}p\r)\ (\mo\ p),\ \sum_{k=0}^{p-1}T_kM_k\eq\f 43\l(\f p3\r)\pmod p.\tag1.8
\endgather$$

\proclaim{Corollary 1.4} Let $p$ be an odd prime. For any integer $x$ we have
$$\sum_{k=0}^{p-1}D_k(x)^2\eq\l(\f{x(x+1)}p\r)\ (\mo\ p).\tag1.9$$
In particular,
$$\sum_{k=0}^{p-1}D_k^2\eq\l(\f2p\r)\ (\mo\ p).\tag1.10$$
\endproclaim
\Proof. It suffices to recall that $D_k(x)=T_k(2x+1,x^2+x)$ and apply Theorem 1.2(ii). \qed
\medskip

\proclaim{Theorem 1.5} Let $b,c\in\Z$ and $d=b^2-4c$.

{\rm (i)} For any $n\in\Z^+$ we have
$$\sum_{k=0}^{n-1}T_k(b,c^2)(b-2c)^{n-1-k}\eq0\ \pmod n\tag1.11$$
and
$$6\sum_{k=0}^{n-1}kT_k(b,c^2)(b-2c)^{n-1-k}\eq0\ \pmod n.\tag1.12$$
If $p$ is an odd prime not dividing $b-2c$, then
$$\f{2c}p\sum_{k=0}^{p-1}\f{T_k(b,c^2)}{(b-2c)^k}\eq-b+(b+2c)\l(\f{b^2-4c^2}p\r)\ (\mo\ p)\tag1.13$$
and
$$\f{12c^2}p\sum_{k=0}^{p-1}\f{kT_k(b,c^2)}{(b-2c)^k}\eq (b+2c)^2\l(1-\l(\f{b^2-4c^2}p\r)\r)-4c^2\ (\mo\ p).\tag1.14$$

{\rm (ii)} Suppose that $d=1$, i.e., there is an $m\in\Z$ such that $b=2m+1$, $c=m^2+m$, and hence $T_k(b,c)=D_k(m)$. Then
$$\f1n\sum_{k=0}^{n-1}(2k+1)T_k(b,c)=\sum_{k=0}^{n-1}\bi n{k+1}\bi{n+k}{k}\l(\f{b-1}2\r)^{k}\in\Z\tag1.15$$
for all $n\in\Z^+$. If $p$ is a prime not dividing $b-1=2m$, then
$$\sum_{k=0}^{p-1}(2k+1)T_k(b,c)\eq p+\f{b+1}{b-1}p\(\l(\f{b+1}2\r)^{p-1}-1\)\ (\mo\ p^3)\tag1.16$$
and
$$\sum_{k=0}^{p-1}(2k+1)^2T_k(b,c)\eq\f2{b-1}\l(\f{(1-b)/2}p\r)=\f1m\l(\f{-m}p\r)\ (\mo\ p).\tag1.17$$
\endproclaim

{\it Example} 1.6.  Putting $b=1$ and $c=\pm1$ in (1.11) we get
$$\sum_{k=0}^{n-1}(-1)^kT_k\eq0\ (\mo\ n)\ \ \t{and}\ \ \sum_{k=0}^{n-1}3^{n-1-k}T_k\eq0\ (\mo\ n),$$
where $n$ is any positive integer.
Also, for a prime $p>3$, (1.13) with $b=1$ and $c=\pm1$ yields $\sum_{k=0}^{p-1}(-1)^kT_k$
and $\sum_{k=0}^{p-1}T_k/3^k$ modulo $p^2$ given by H. Q. Cao and H. Pan [CP].

\medskip

\Remark\ 1.7. For any $n\in\Z^+$, we have
$$\f1n\sum_{k=0}^{n-1}(2k+1)T_k3^{n-1-k}=\sum_{k=0}^{n-1}\bi{n-1}k(-1)^{n-1-k}(k+1)\bi{2k}k,$$
for, if $a_n$ denotes the left-hand side or the right-hand side of the last equality, then by the Zeilberger algorithm
[PWZ, pp.\,101-119], we have the recurrence
$$(n+1)(2n+1)a_{n+2}=(4n^2+10n+3)a_{n+1}+3n(2n+3)a_n,\quad n=0,1,2,\ldots.$$
If $b,c\in\Z$ with $b^2-4c=1$, then for any prime $p\nmid c$, by (1.16) we have
$$\sum_{k=0}^{p-1}(2k+1)T_k(b,c)\eq p\ (\mo\ p^2).$$
\medskip

\proclaim{Theorem 1.8} Let $b,c\in\Z$ and $d=b^2-4c$.

{\rm (i)} For any $n\in\Z^+$ we have
$$\sum_{k=0}^{n-1}(2k+1)T_k(b,c)^2(-d)^{n-1-k}\eq0\ \ (\mo\ n),\tag1.18$$
and furthermore
$$b\sum_{k=0}^{n-1}(2k+1)T_k(b,c)^2(-d)^{n-1-k}=nT_n(b,c)T_{n-1}(b,c).\tag1.19$$

{\rm (ii)}  For any $n\in\Z^+$ we have
$$\f1{n^2}\sum_{k=0}^{n-1}(2k+1)T_k(b,c)^2d^{n-1-k}=\sum_{k=0}^{n-1}\bi{n-1}k\bi{n+k}kC_kc^kd^{n-1-k}\in\Z.\tag1.20$$
If $c$ is nonzero and $p$ is an odd prime not dividing $d$, then
$$\f1{p^2}\sum_{k=0}^{p-1}(2k+1)\f{T_k(b,c)^2}{d^k}\eq1+\f{b^2}c\cdot\f{(\f{d}p)-1}2\ (\mo\ p).\tag1.21$$
\endproclaim

Now we give one more theorem.

\proclaim{Theorem 1.9} Let $p>3$ be a prime. Then
$$\align\sum_{k=0}^{p-1}\f{T_k(6,-3)^2}{48^k}\eq&\l(\f{-1}p\r)+\f{p^2}3E_{p-3}\ (\mo\ p^3),\tag1.22
\\ \sum_{k=0}^{p-1}\f{T_k(2,-1)^2}{8^k}\eq&\l(\f{-2}p\r)\ (\mo\ p^2),\tag1.23
\\ \sum_{k=0}^{p-1}\f{T_k(2,-3)^2}{16^k}\eq&\l(\f{p}3\r)\ (\mo\ p^2),\tag1.24
\\\sum_{k=1}^{p-1}\f{D_k^2}{k^2}\eq&-2q_p(2)^2\pmod{p},\tag1.25
\endalign$$
where $E_0,E_1,E_2,\ldots$ are Euler numbers, and $q_p(2)$ denotes the Fermat quotient $(2^{p-1}-1)/p$.
\endproclaim
\Remark\ 1.10. (1.25) was conjectured by the author in [Su3].
\medskip

We will show Theorems 1.2 and 1.5 in Sections 2 and 3 respectively. Section 4 is devoted to our proofs of Theorems 1.8 and 1.9.
In Section 5 we are going to pose more conjectures for further research.

\heading{2. Proof of Theorem 1.2}\endheading

The following lemma essentially follows from [ST, (1.5)], but we will give a direct proof.
\proclaim{Lemma 2.1} Let $p$ be an odd prime and let $m\in\Z$ with
$m\not\eq0\ (\mo\ p)$. Then
$$\sum_{k=0}^{(p-1)/2}\f{\bi{2k}k}{m^k}\eq\l(\f {m(m-4)}p\r)\ (\mo\ p)\tag2.1$$
and
$$\sum_{k=0}^{(p-1)/2}\f{C_k}{m^k}\eq\f m2-\f{m-4}2\l(\f {m(m-4)}p\r)\ (\mo\ p).\tag2.2$$
\endproclaim
\Proof. Clearly
$$\bi{2k}k=\bi{-1/2}k(-4)^k\eq\bi{(p-1)/2}k(-4)^k\ (\mo\ p)$$
for all $k=0,\ldots,p-1$. Thus
$$\align\sum_{k=0}^{(p-1)/2}\f{\bi{2k}k}{m^k}\eq&\sum_{k=0}^{(p-1)/2}\bi{(p-1)/2}k\f{(-4)^k}{m^k}
=\l(1-\f 4m\r)^{(p-1)/2}
\\=&\f{(m(m-4))^{(p-1)/2}}{m^{p-1}}\eq\l(\f {m(m-4)}p\r)\ (\mo\ p).
\endalign$$
This proves (2.1).

Observe that
$$\align \sum_{k=0}^{(p-1)/2}\f{\bi{2k+1}k}{m^k}
=&\f{\bi
p{(p-1)/2}}{m^{(p-1)/2}}+\f12\sum_{k=0}^{(p-3)/2}\f{\bi{2k+2}{k+1}}{m^k}
\\\eq&\f m2\sum_{k=0}^{(p-1)/2}\f{\bi{2k}k}{m^k}-\f m2\ (\mo\ p).
\endalign$$
Hence
$$\align\sum_{k=0}^{(p-1)/2}\f{C_k}{m^k}=&\sum_{k=0}^{(p-1)/2}\f{2\bi{2k}k-\bi{2k+1}{k}}{m^k}
\\\eq& \l(2-\f m2\r)\sum_{k=0}^{(p-1)/2}\f{\bi{2k}k}{m^k}+\f m2
\\\eq&\f m2-\f{m-4}2\l(\f{m(m-4)}p\r)\ (\mo\ p).\endalign$$
So (2.2) also holds. We are done. \qed

\medskip
\noindent{\it Proof of Theorem} 1.2(i). In the case $c\eq0\ (\mo\
p)$, as $T_k(b,c)\eq b^k\ (\mo\ c)$ for all $k\in\N$, we have
$$\sum_{k=0}^{p-1}\f{T_k(b,c)}{m^k}\eq\sum_{k=0}^{p-1}\f {b^k}{m^k}\eq\l(\f{(m-b)^2}p\r)\ (\mo\ p).$$
So (1.1) holds if $p\mid c$. Note that (1.2) is trivial when $p\mid
c$.

Suppose that $c\not\eq0\ (\mo\ p)$. For any $n\in\N$, clearly,
$$T_n(b,c)=\sum_{k=0}^{\lfloor n/2\rfloor}\bi n{2k}\bi{2k}kb^{n-2k}c^k
\eq\cases \bi n{n/2}c^{n/2}\ (\mo\ b)&\t{if}\ 2\mid n,\\0\ (\mo\
b)&\t{if}\ 2\nmid n,\endcases$$
and similarly,
$$M_n(b,c)=\sum_{k=0}^{\lfloor n/2\rfloor}\bi n{2k}C_kb^{n-2k}c^k
\eq\cases C_{n/2}c^{n/2}\ (\mo\ b)&\t{if}\ 2\mid n,\\0\ (\mo\
b)&\t{if}\ 2\nmid n.\endcases$$
In the case $b\eq0\ (\mo\ p)$, by
applying Lemma 2.1 we obtain
$$\sum_{k=0}^{p-1}\f{T_k(b,c)}{m^k}\eq\sum_{k=0}^{(p-1)/2}\f{\bi{2k}kc^k}{m^{2k}}
\eq\sum_{k=0}^{(p-1)/2}\f{\bi{2k}k}{(m^2c^{p-2})^k}\eq\l(\f{m^2-4c}p\r)\
(\mo\ p)$$ and
$$\align\sum_{k=0}^{p-1}\f{M_k(b,c)}{m^k}\eq&\sum_{k=0}^{(p-1)/2}\f{C_kc^k}{m^{2k}}
\eq\sum_{k=0}^{(p-1)/2}\f{C_k}{(m^2c^{p-2})^k}
\\\eq&\f{m^2}{2c}-\f{m^2-4c}{2c}\l(\f{m^2-4c}p\r)\ (\mo\ p).
\endalign$$
So (1.1) and (1.2) hold when $p\mid b$.

Below we assume that $p\nmid bc$. Observe that
$$\align \sum_{n=0}^{p-1}\f{T_n(b,c)}{m^n}=&\sum_{n=0}^{p-1}\f1{m^n}\sum_{k=0}^{\lfloor n/2\rfloor}\bi n{2k}\bi{2k}kb^{n-2k}c^k
\\=&\sum_{k=0}^{(p-1)/2}\bi{2k}k\f{c^k}{b^{2k}}\sum_{n=0}^{p-1}\f{b^n}{m^n}\bi n{2k};
\endalign$$
similarly,
$$\sum_{n=0}^{p-1}\f{M_n(b,c)}{m^n}=\sum_{k=0}^{(p-1)/2}C_k\f{c^k}{b^{2k}}\sum_{n=0}^{p-1}\f{b^n}{m^n}\bi n{2k}.$$

Now we consider the case $m\eq b\ (\mo\ p)$. For
$k\in\{0,1,\ldots,(p-1)/2\}$ we have
$$\sum_{n=0}^{p-1}\f{b^n}{m^n}\bi n{2k}\eq\sum_{n=2k}^{p-1}\bi n{2k}=\bi p{2k+1}\ \ (\mo\ p)$$
with the help of a well-known identity of Chu (see, (1.52) of H. Gould [G, p,\,7] or (5.26) of [GKP, p.\,169]). Thus, by the above,
$$\sum_{n=0}^{p-1}\f{T_n(b,c)}{m^n}\eq\bi{p-1}{(p-1)/2}\f{c^{(p-1)/2}}{b^{p-1}}\eq\l(\f{-c}p\r)=\l(\f{(m-b)^2-4c}p\r)\ (\mo\ p)$$
and
$$\sum_{n=0}^{p-1}\f{M_n(b,c)}{m^n}\eq C_{(p-1)/2}\f{c^{(p-1)/2}}{b^{p-1}}
\eq2\l(\f{-c}p\r)=2\l(\f{(m-b)^2-4c}p\r)\ (\mo\ p).$$ So (1.1) and
(1.2) are true.

Below we consider the remaining case $m\not\eq b\ (\mo\ p)$. Observe
that
$$\align&\sum_{n=0}^{p-1}\f{b^n}{m^n}\bi n{2k}=[x^{2k}]\sum_{n=0}^{p-1}\f{b^n}{m^n}(1+x)^n
\\\eq&[x^{2k}]\sum_{n=0}^{p-1}(b+bx)^nm^{p-1-n}=[x^{2k}]\f{(b+bx)^p-m^p}{b+bx-m}
\\=&[x^{2k}]\f{(b+bx)^p-m^p}{-(m-b)^p}\cdot\f{(bx)^p-(m-b)^p}{bx-(m-b)}
\\\eq&[x^{2k}]\f{b^p+b^px^p-m^p}{-(m-b)^p}\sum_{j=0}^{p-1}(bx)^j(m-b)^{p-1-j}\eq \f{b^{2k}}{(m-b)^{2k}}\ (\mo\ p).
\endalign$$
Therefore, with the help of Lemma 2.1,
$$\sum_{k=0}^{p-1}\f{T_n(b,c)}{m^n}
\eq\sum_{k=0}^{(p-1)/2}\bi{2k}k\f{c^k}{b^{2k}}\cdot\f{b^{2k}}{(m-b)^{2k}}\eq\l(\f{(m-b)^2-4c}p\r)\ (\mo\ p).
$$
This proves (1.1).

In a similar way,
$$\sum_{n=0}^{p-1}\f{M_n(b,c)}{m^n}\eq\sum_{k=0}^{(p-1)/2}C_k\f{c^k}{(m-b)^{2k}}
\eq\sum_{k=0}^{(p-1)/2}\f{C_k}{M^k}\ (\mo\ p),$$ where
$M:=(m-b)^2c^{p-2}$. Applying Lemma 2.1 we get the desired (1.2).
\qed

\proclaim{Lemma 2.2} Let $b,c\in\Z$ and $d=b^2-4c$. Let $p$ be any
odd prime and let $n\in\{0,\ldots,p-1\}$. If $p\nmid d$ or
$p/2<n<p$, then
$$T_n(b,c)\eq\l(\f dp\r)d^nT_{p-1-n}(b,c)\ \ (\mo\ p).\tag2.3$$
\endproclaim
\Proof. If $p\mid d$, then
$$T_n(b,c)\eq [x^n]\l(x^2+bx+\f{b^2}4\r)^n=[x^n]\l(x+\f b2\r)^{2n}=\bi{2n}n\f{b^n}{2^n}\ (\mo\ p).$$
Note that for $n=(p+1)/2,\ldots,p-1$ we have
$$\bi{2n}n=\f{(2n)!}{(n!)^2}\eq0\ (\mo\ p).$$

Now assume that $p\nmid d$. Then
$$\align &d^nT_{p-1-n}(b,c)=d^n(\sqrt d)^{p-1-n}P_{p-1-n}\l(\f b{\sqrt d}\r)
\\=&d^{(p-1)/2}\sum_{k=0}^{p-1-n}\bi{p-1-n+k}{2k}\bi{2k}k\l(\f{b/\sqrt d-1}2\r)^k(\sqrt d)^n
\\=&d^{(p-1)/2}\sum_{k=0}^{p-1}\bi{n+k-p}{2k}\bi{2k}k\l(\f{b-\sqrt d}{2\sqrt d}\r)^k(\sqrt d)^n
\\\eq&d^{(p-1)/2}\sum_{k=0}^{n}\bi{n+k}{2k}\bi{2k}k\l(\f{b-\sqrt d}{2\sqrt d}\r)^k(\sqrt d)^n
\\\eq&\l(\f dp\r)(\sqrt d)^nP_n\l(\f b{\sqrt d}\r)
=\l(\f dp\r)T_n(b,c)\ \ (\mo\ p).
\endalign$$
This concludes the proof. \qed

\Remark\ 2.3. Lemma 2.2 in the case $p\nmid d$ is essentially known
(see, e.g., [N, (14)]), but our proof is simple and direct. By Lemma
2.2, for any prime $p>3$ we have
$$\sum_{k=0}^{p-1}\f{T_k^2}{9^k}=\sum_{k=0}^{p-1}\l(\f{T_k}{(-3)^k}\r)^2
\eq\sum_{k=0}^{p-1}\l(\l(\f{-3}p\r)T_{p-1-k}\r)^2=\sum_{j=0}^{p-1}T_j^2\
\ (\mo\ p)$$
and hence $\sum_{k=0}^{p-1}T_k^2/9^k\eq(\f{-1}p)\ (\mo\ p)$ in light of Example 1.3.
\medskip

Let $A$ and $B$ be integers. The Lucas sequence $u_n=u_n(A,B)\
(n\in\N)$ is defined by
$$u_0=0,\ u_1=1,\ \t{and}\ u_{n+1}=Au_n-Bu_{n-1}\ (n=1,2,3,\ldots).$$
Let $\al$ and $\beta$ be the two roots of the equation $x^2-Ax+B=0$.
It is well-known that if $\Delta=A^2-4B\not=0$ then
$$u_n=\f{\al^n-\beta^n}{\al-\beta}\quad\ \t{for all}\ n=0,1,2,\ldots.$$

\proclaim{Lemma 2.4} Let $A$ and $B$ be integers. For any odd prime
$p$ we have
$$u_p(A,B)\eq\l(\f{A^2-4B}p\r)\ (\mo\ p).$$
\endproclaim
\Proof. Though this is a known result, here we provide a simple
proof.

If $\Delta=A^2-4B\eq0\ (\mo\ p)$, then
$$u_n(A,B)\eq u_n\l(A,\f{A^2}4\r)=n\l(\f A2\r)^{n-1}\ (\mo\ p)\ \ \ \t{for}\ n=1,2,3,\ldots$$
and in particular $u_p(A,B)\eq0\ (\mo\ p)$.

 When $\Delta\not\eq0\ (\mo\ p)$, we have
$$\Delta u_p(A,B)=(\al-\beta)(\al^p-\beta^p)\eq(\al-\beta)(\al-\beta)^p=\Delta^{(p+1)/2}\ (\mo\ p)$$
with $\al$ and $\beta$ the two roots of the equation $x^2-Ax+B=0$,
hence $u_p(A,B)\eq(\f{\Delta}p)\ (\mo\ p)$ as desired. \qed

\medskip
\noindent{\it Proof of Theorem} 1.2(ii). Suppose that
$d=b^2-4c\not\eq0\pmod p$. By Lemma 2.2,
$$\align\l(\f dp\r)\sum_{k=0}^{p-1}\f{T_k(b,c)^2}{d^k}
\eq&\sum_{k=0}^{p-1}T_k(b,c)T_{p-1-k}(b,c)=[x^{p-1}]\(\sum_{n=0}^\infty T_n(b,c)x^n\)^2
\\=&[x^{p-1}]\f1{1-2bx+dx^2}=[x^p]\f x{1-2bx+dx^2}\ \ (\mo\ p).
\endalign$$
Write
$$\f x{1-2bx+dx^2}=\sum_{n=0}^\infty u_nx^n.$$
Then $u_0=0$ and $u_1=1$. Since $(1-2bx+dx^2)\sum_{n=0}^\infty
u_nx^n=x$, we have $u_n-2bu_{n-1}+du_{n-2}=0$ for $n=2,3,\ldots$,
hence $u_n=u_n(2b,d)$ for all $n\in\N$. Thus, with the help of Lemma
2.3, from the above we obtain
$$\l(\f dp\r)\sum_{k=0}^{p-1}\f{T_k(b,c)^2}{d^k}\eq u_p(2b,d)\eq\l(\f{4b^2-4d}p\r)=\l(\f cp\r)\ (\mo\ p).$$
This proves (1.3).

Now suppose that $b\not\eq 2c\ (\mo\ p)$ and set
$D=b^2-4c^2=(b-2c)(b+2c)$. If $p\mid D$, then $b\eq-2c\not\eq0\
(\mo\ p)$ and $T_k(b,c^2)\eq[x^k](x^2+bx+b^2/4)^k=[x^k](x+b/2)^{2k}$, hence
$$\sum_{k=0}^{p-1}\f{T_k(b,c^2)^2}{(b-2c)^{2k}}\eq\sum_{k=0}^{p-1}\f{(\bi{2k}k(b/2)^k)^2}{(2b)^{2k}}=\sum_{k=0}^{p-1}\f{\bi{2k}k^2}{16^k}
\eq\l(\f{-1}p\r)\ (\mo\ p).$$ The last step can be easily explained
as follows:
$$\align\sum_{k=0}^{p-1}\f{\bi{2k}k^2}{16^k}\eq&\sum_{k=0}^{(p-1)/2}\bi{-1/2}k^2
\\\eq&\sum_{k=0}^{(p-1)/2}\bi{(p-1)/2}k\bi{(p-1)/2}{(p-1)/2-k}
\\=&[x^{(p-1)/2}](1+x)^{(p-1)/2+(p-1)/2}
\\=&\bi{p-1}{(p-1)/2}\eq\l(\f{-1}p\r)\ (\mo\ p).
\endalign$$
Below we assume that $p\nmid D$. By Lemma 2.2 and Fermat's little
theorem,
 $$\l(\f Dp\r)\sum_{k=0}^{p-1}\f{T_k(b,c^2)^2}{(b-2c)^{2k}}\eq C\ (\mo\ p),$$
 where
 $$\align C=&\sum_{k=0}^{p-1}D^kT_k(b,c^2)(b-2c)^{2(p-1-k)}T_{p-1-k}(b,c^2)
\\ =&[x^{p-1}]\(\sum_{k=0}^\infty T_k(b,c^2)(Dx)^k\)\sum_{l=0}^\infty T_l(b,c^2)(b-2c)^{2l}x^l
\\=&[x^{p-1}]\f1{\sqrt{1-2b(Dx)+D(Dx)^2}}
\cdot\f1{\sqrt{1-2b(b-2c)^2x+D(b-2c)^4x^2}}
\\=& [y^{p-1}]\f{(b-2c)^{p-1}}{\sqrt{(1-2b(b+2c)y+(b+2c)^2Dy^2)(1-2b(b-2c)y+D(b-2c)^2y^2)}}.
\endalign$$
(Note that $y$ corresponds to $(b-2c)x$.) Therefore
$$\align C\eq&[y^{p-1}]\f1{1-Dy}\cdot\f1{\sqrt{(1-(b+2c)^2y)(1-(b-2c)^2y)}}
\\=&[y^{p-1}]\sum_{n=0}^\infty(Dy)^n\f1{\sqrt{1-2(b^2+4c^2)y+D^2y^2}}\ \ (\mo\ p).
\endalign$$
Observe that $(b^2+4c^2)^2-4(4b^2c^2)=(b^2-4c^2)^2=D^2$ and hence
$$\f1{\sqrt{1-2(b^2+4c^2)y+D^2y^2}}=\sum_{k=0}^\infty T_k(b^2+4c^2,4b^2c^2)y^k.$$
So we have
$$\align C\eq&\sum_{k=0}^{p-1}T_k(b^2+4c^2,4b^2c^2)D^{p-1-k}\eq\sum_{k=0}^{p-1}\f{T_k(b^2+4c^2,4b^2c^2)}{D^k}
\\\eq&\l(\f{(D-(b^2+4c^2))^2-4(4b^2c^2)}p\r)=\l(\f{-16c^2D}p\r)\ (\mo\ p)
\endalign$$
with the help of the first part of Theorem 1.2.

Combining the above, we finally obtain (1.4). We are done. \qed

\proclaim{Lemma 2.5} Let $b$ and $c$ be integers. For any odd prime
$p$, we have
$$T_p(b,c)\eq b\ (\mo\ p),\ \ T_{p+1}(b,c)\eq b^2\ (\mo\ p),\tag2.4$$
and
$$T_{p-1}(b,c)\eq\l(\f{b^2-4c}p\r)\ (\mo\ p).\tag2.5$$
\endproclaim
\Proof. Since $\bi pk\eq0\ (\mo\ p)$ for all $k=1,\ldots,p-1$, we
have
$$T_p(b,c)=\sum_{k=0}^{(p-1)/2}\bi p{2k}\bi{2k}kb^{p-2k}c^k\eq\bi p0 b^p\eq b\ (\mo\ p)$$
with the help of Fermat's little theorem. If $1<k<p$, then
$$\bi{p+1}k=\f{p(p+1)}{k(k-1)}\bi{p-1}{k-2}\eq0\ \ (\mo\ p).$$
Thus
$$\align T_{p+1}(b,c)=&\sum_{k=0}^{(p+1)/2}\bi{p+1}k\bi{p+1-k}kb^{p+1-2k}c^k
\\\eq& b^{p+1}+\bi{p+1}1\bi{p}1 b^{p-1}c\eq b^2\ (\mo\ p).
\endalign$$

If $p\mid b$, then (2.5) is valid since
$$\align T_{p-1}(b,c)=&\sum_{k=0}^{(p-1)/2}\bi{p-1}{2k}\bi{2k}kb^{p-1-2k}c^k
\\\eq&\bi{p-1}{(p-1)/2}c^{(p-1)/2}\eq\l(\f{-c}p\r)=\l(\f{b^2-4c}p\r)\ (\mo\ p).
\endalign$$
When $p\nmid b$, we have
$$\align T_{p-1}(b,c)\eq&\sum_{k=0}^{(p-1)/2}\bi{2k}k\f{c^k}{b^{2k}}=\sum_{k=0}^{(p-1)/2}\bi{-1/2}k(-4)^k\f{c^k}{b^{2k}}
\\\eq&\sum_{k=0}^{(p-1)/2}\bi{(p-1)/2}k\l(-\f{4c}{b^2}\r)^k=\l(1-\f{4c}{b^2}\r)^{(p-1)/2}
\\\eq&\l(\f{b^2-4c}p\r)\ (\mo\ p).
\endalign$$
This concludes the proof. \qed

\medskip
\noindent{\it Proof of Theorem} 1.2(iii). Suppose that
$d=b^2-4c\not\eq0\pmod p$. By Lemma 2.2,
$$\sum_{k=0}^{p-1}\f{T_k(b,c)M_k(b,c)}{d^k}\eq\l(\f dp\r)S_1\ (\mo\ p)$$
where
$$\align S_1=&\sum_{k=0}^{p-1}T_{p-1-k}(b,c)M_k(b,c)=[x^{p-1}]\sum_{j=0}^\infty T_j(b,c)x^j\sum_{k=0}^\infty M_k(b,c)x^k
\\=&[x^{p-1}]\f1{\sqrt{1-2bx+dx^2}}\times\f{1-bx-\sqrt{1-2bx+dx^2}}{2cx^2}
\\=&\f1{2c}[x^{p+1}]\l(\f{1-bx}{\sqrt{1-2bx+dx^2}}-1\r)=\f{T_{p+1}(b,c)-bT_p(b,c)}{2c}.
\endalign$$
In light of Lemma 2.5,  $S_1\eq0\ (\mo\ p)$ and hence (1.5) follows.

Now suppose that $D=b^2-4c^2\not\eq0\pmod p$.  In view of Lemma 2.2
and Fermat's little theorem,
$$\align&\sum_{k=0}^{p-1}\f{T_k(b,c^2)M_k(b,c^2)}{(b-2c)^{2k}}
\\\eq&\l(\f Dp\r)\sum_{k=0}^{p-1}\f{D^kT_{p-1-k}(b,c^2)}{(b-2c)^{2k}}M_k(b,c^2)
\eq \l(\f Dp\r)S_2\ (\mo\ p),
\endalign$$
where
$$\align S_2=&\sum_{k=0}^{p-1}(b-2c)^{p-1-k}T_{p-1-k}(b,c^2)M_k(b,c^2)(b+2c)^k
\\=&[x^{p-1}]\sum_{j=0}^\infty T_j(b,c^2)((b-2c)x)^j\sum_{k=0}^\infty M_k(b,c^2)((b+2c)x)^k
\\=&[x^{p-1}]\f{1-b(b+2c)x-\sqrt{1-2b(b+2c)x+D(b+2c)^2x^2}}{2c^2((b+2c)x)^2\sqrt{1-2b(b-2c)x+D(b-2c)^2x^2}}
\\=&\f1{2c^2(b+2c)^2}[x^{p+1}]\f{1-b(b+2c)x}{\sqrt{1-2b(b-2c)x+D(b-2c)^2x^2}}
\\&-\f1{2c^2(b+2c)^2}[x^{p+1}]\f{\sqrt{(1-Dx)(1-(b+2c)^2x)}}{\sqrt{(1-Dx)(1-(b-2c)^2x)}}.
\endalign$$
Recall the identity $(b^2+4c^2)^2-4(4b^2c^2)=D^2$ and observe that
$$\align 2c^2(b+2c)^2S_2=&[y^{p+1}]\f{(b-2c)^{p+1}}{\sqrt{1-2by+Dy^2}}-b(b+2c)[y^p]\f{(b-2c)^p}{\sqrt{1-2by+Dy^2}}
\\&-[x^{p+1}]\f{1-(b+2c)^2x}{\sqrt{1-2(b^2+4c^2)x+D^2x^2}}
\\\eq&(b-2c)^2T_{p+1}(b,c^2)-b(b+2c)(b-2c)T_p(b,c^2)
\\&-T_{p+1}(b^2+4c^2,4b^2c^2)+(b+2c)^2T_p(b^2+4c^2,4b^2c^2)\ (\mo\ p).
\endalign$$
Applying Lemma 2.5 we get
$$\align 2c^2(b+2c)^2S_2\eq&(b-2c)^2b^2-b^2D-(b^2+4c^2)^2+(b+2c)^2(b^2+4c^2)
\\=&8bc^2(b+2c)\ \ (\mo\ p).
\endalign$$
Thus $S_2\eq 4b/(b+2c)\ (\mo\ p)$ and this concludes the proof of
(1.6). \qed

\heading{3. Proof of Theorem 1.5}\endheading

\proclaim{Lemma 3.1} Let $b$ and $c$ be integers. For all
$n=1,2,3,\ldots$ we have
$$2c\sum_{k=0}^{n-1}T_k(b,c^2)(b-2c)^{n-1-k}
=-nT_n(b,c^2)+(b+2c)nT_{n-1}(b,c^2).\tag3.1$$
\endproclaim
\Proof. In the case $n=1$ both sides of (3.1) coincide with $2c$.
Denote by $f(n)$ the right-hand side of (3.1). Clearly it suffices
to show that for any positive integer $n$ we have
$$\align &f(n+1)-(b-2c)f(n)
\\=&2c\sum_{k=0}^{n}T_k(b,c^2)(b-2c)^{n-k}-2c\sum_{k=0}^{n-1}T_k(b,c^2)(b-2c)^{n-k}=2cT_n(b,c^2).
\endalign$$
Observe that
$$\align &f(n+1)-(b-2c)f(n)
\\=&-(n+1)T_{n+1}(b,c^2)+(b+2c)(n+1)T_n(b,c^2)
\\&-(b-2c)\l(-nT_n(b,c^2)+(b+2c)nT_{n-1}(b,c^2)\r)
\\=&-(n+1)T_{n+1}(b,c^2)+(4c^2-b^2)nT_{n-1}(b,c^2)
\\&+\l(n(b-2c)+(n+1)(b+2c)\r)T_n(b,c^2)
\\=&-(2n+1)bT_n(b,c^2)+\l(n(b-2c)+(n+1)(b+2c)\r)T_n(b,c^2)=2cT_n(b,c^2)
\endalign$$
with the help of the recursion for $T_n(b,c^2)$.

The above proof of (3.1) is simple. However, the reader might wonder
how (3.1) was found. Set $D=b^2-4c^2$. Then
$$\align\sum_{k=0}^{n-1}T_k(b,c^2)(b-2c)^{n-1-k}
=&[x^{n-1}]\f1{\sqrt{1-2bx+Dx^2}}\cdot\f1{1-(b-2c)x}
\\=&[x^{n-1}](1-(b-2c)x)^{-3/2}(1-(b+2c)x)^{-1/2}
\endalign$$
 and hence
 $$-2c\sum_{k=0}^{n-1}T_k(b,c^2)(b-2c)^{n-1-k}=[x^{n-1}]\f d{dx}\sqrt{\f{1-(b+2c)x}{1-(b-2c)x}}.$$
Observe that
$$\align \sqrt{\f{1-(b+2c)x}{1-(b-2c)x}}=&\f{1-(b+2c)x}{\sqrt{1-2bx+Dx^2}}=(1-(b+2c)x)\sum_{k=0}^\infty T_k(b,c^2)x^k
\\=&1+\sum_{k=1}^\infty(T_k(b,c^2)-(b+2c)T_{k-1}(b,c^2))x^k
\endalign$$ and thus
$$[x^{n-1}]\f d{dx}\sqrt{\f{1-(b+2c)x}{1-(b-2c)x}}=n\l(T_n(b,c^2)-(b+2c)T_{n-1}(b,c^2)\r).$$
Therefore (3.1) follows. \qed

\proclaim{Lemma 3.2} Let $b\in\Z$, $c\in\Z\sm\{0\}$ and $n\in\Z^+$. Then
$$\aligned&\f 3n\sum_{k=0}^{n-1}kT_k(b,c^2)(b-2c)^{n-1-k}-\sum_{k=0}^{n-1}T_k(b,c^2)(b-2c)^{n-1-k}
\\&\qquad=\f{(b+4c)T_n(b,c^2)-(b+2c)^2T_{n-1}(b,c^2)}{4c^2}.
\endaligned\tag3.2$$
\endproclaim
\Proof. Note that for any $k\in\N$ we have
$$T_k(2c,c^2)=[x^k](x^2+2cx+c^2)^k=[x^k](x+c)^{2k}=\bi{2k}kc^k.$$
In the case $b=2c$, we can easily verify that both sides of (3.2) coincide with $(2-3/n)\bi{2n-2}{n-1}c^{n-1}$.

Below we assume $b\not=2c$ and define
$$\sigma_n:=\sum_{k=0}^{n-1}(n-k)T_k(b,c^2)(b-2c)^{n-1-k}.$$
Clearly
$$\sigma_n=[x^{n-1}]\(\sum_{k=0}^\infty T_k(b,c^2)x^k\)\sum_{l=0}^\infty(l+1)(b-2c)^lx^l.$$
For $|z|<1$ we have
$$\f1{(1-z)^2}=\sum_{l=0}^\infty\bi{-2}l(-z)^l=\sum_{l=0}^\infty\bi{l+1}lz^l.$$
Thus
$$\align\sigma_n=&[x^{n-1}]\f1{\sqrt{1-2bx+(b^2-4c^2)x^2}}\times\f1{(1-(b-2c)x)^2}
\\=&[x^{n-1}](1-(b+2c)x)^{-1/2}(1-(b-2c)x)^{-5/2}=[x^{n-1}]\f d{dx}f(x),
\endalign$$
where
$$\align f(x)=&\(-\f{b(b+2c)}{12c^2(b-2c)}+\f{(b+2c)^2}{12c^2}x+\f2{3(b-2c)(1-(b-2c)x)}\)
\\&\times\f1{\sqrt{1-2bx+(b^2-4c^2)x^2}}
\\=&\(-\f{b(b+2c)}{12c^2(b-2c)}+\f{(b+2c)^2}{12c^2}x+\f2{3(b-2c)}\sum_{j=0}^\infty(b-2c)^jx^j\)
\\&\times\sum_{k=0}^\infty T_k(b,c^2)x^k.
\endalign$$
Therefore
$$\align\f{\sigma_n}n=[x^n]f(x)=&-\f{b(b+2c)}{12c^2(b-2c)}T_n(b,c^2)+\f{(b+2c)^2}{12c^2}T_{n-1}(b,c^2)
\\&+\f2{3(b-2c)}\sum_{k=0}^n T_k(b,c^2)(b-2c)^{n-k},
\endalign$$
i.e.,
$$\align &\sum_{k=0}^{n-1}T_k(b,c^2)(b-2c)^{n-1-k}-\f1n\sum_{k=0}^{n-1}kT_k(b,c^2)(b-2c)^{n-1-k}
\\=&\f23\sum_{k=0}^{n-1}T_k(b,c^2)(b-2c)^{n-1-k}+\f23\cdot\f{T_n(b,c^2)}{b-2c}
\\&+\f{b+2c}{12c^2(b-2c)}\l((b^2-4c^2)T_{n-1}(b,c^2)-bT_n(b,c^2)\r).
\endalign$$
This yields the desired (3.2). \qed

\medskip
\noindent{\it Proof of Theorem} 1.5(i). Let $n$ be any positive
integer. Since $T_k(b,0)=[x^k]x^k(x+b)^k=b^k$ for all $k\in\N$,
(1.11) and (1.12) hold trivially when $c=0$.

Now assume that $c\not=0$. By Lemma 3.1 we have
$$\f1n\sum_{k=0}^{n-1}T_k(b,c^2)(b-2c)^{n-1-k}
=\f{bT_{n-1}(b,c^2)-T_n(b,c^2)}{2c}+T_{n-1}(b,c^2).$$ Observe
that
$$\align &T_n(b,c^2)-bT_{n-1}(b,c^2)
\\=&\sum_{k\in\N}\bi n{2k}\bi{2k}kb^{n-2k}(c^2)^k-\sum_{k\in\N}\bi {n-1}{2k}\bi{2k}kb^{n-2k}(c^2)^k
\\=&\sum_{k=1}^n\bi{n-1}{2k-1}\bi{2k}kb^{n-2k}c^{2k}
=2c\sum_{k=1}^n\bi{n-1}{2k-1}\bi{2k-1}{k-1}b^{n-2k}c^{2k-1}
\\=&2c\sum_{0<k\ls\lfloor n/2\rfloor}\bi{n-1}{k-1}\bi{n-k}kb^{n-2k}c^{2k-1}\eq0\pmod{2c}.
\endalign$$
Therefore (1.11) holds.
In light of Lemma 3.2, (1.12) is reduced to the congruence
$$(b+4c)T_n(b,c^2)\eq (b+2c)^2T_{n-1}(b,c^2)\pmod{2c^2}.$$
In fact, as $\bi{2k}k=2\bi{2k-1}{k-1}$ for all $k\in\Z^+$, we have
$$\align &(b+4c)T_n(b,c^2)-(b+2c)^2T_{n-1}(b,c^2)
\\=&(b+4c)\sum_{k=0}^{\lfloor n/2\rfloor}\bi n{2k}\bi{2k}kb^{n-2k}c^{2k}
\\&-(b+2c)^2\sum_{k=0}^{\lfloor(n-1)/2\rfloor}\bi{n-1}{2k}\bi{2k}kb^{n-1-2k}c^{2k}
\\\eq&(b+4c)b^n-(b+2c)^2b^{n-1}\eq0\pmod{2c^2}.
\endalign$$
So (1.12) is valid.

Now write $D=b^2-4c^2$ and suppose that $p$ is an odd prime not dividing $b-2c$.
In view of Lemmas 2.4 and 3.1 and Fermat's little theorem, we have
$$\align\f{2c}p\sum_{k=0}^{p-1}\f{T_k(b,c^2)}{(b-2c)^k}=&\f{(b+2c)T_{p-1}(b,c^2)-T_p(b,c^2)}{(b-2c)^{p-1}}
\\\eq&(b+2c)\l(\f{D}p\r)-b\ (\mo\ p).
\endalign$$
This proves (1.13). If $p\mid c$, then $(\f Dp)=(\f{b^2}p)=1$ and hence (1.14) becomes obvious.
When $p\nmid c$, by (3.2), (1.11) and Lemma 2.5 we get
$$\align \f3p\sum_{k=0}^{p-1}\f{kT_k(b,c^2)}{(b-2c)^k}\eq&\f{(b+4c)T_p(b,c^2)-(b+2c)^2T_{p-1}(b,c^2)}{4c^2}
\\\eq&\f{(b+4c)b-(b+2c)^2(\f Dp)}{4c^2}\pmod{p}
\endalign$$
and hence (1.14) follows.
\qed

\proclaim{Lemma 3.3} For $k\in\N$ and $n\in\Z^+$ we have
$$\sum_{m=0}^{n-1}(2m+1)^2\bi{m+k}{2k}=(4n^2-1)\f{n-k}{2k+3}\bi{n+k}{2k}.\tag3.3$$
\endproclaim
\Proof. Observe that
$$\align&(4n^2-1)\f{n-k}{2k+3}\bi{n+k}{2k}+(2n+1)^2\bi{n+k}{2k}
\\=&(4n^2+8n+3)\f{n+1+k}{2k+3}\bi{n+k}{2k}
\\=&(4(n+1)^2-1)\f{n+1-k}{2k+3}\bi{n+1+k}{2k}.
\endalign$$
So we can easily prove (3.3) by induction on $n$. \qed

\medskip
\noindent{\it Proof of Theorem} 1.5(ii). We prove (1.15) by induction. (1.15) is obvious when
$n=1$.

Now suppose the validity of (1.15) for a fixed $n\in\Z^+$. Observe
that
$$\align&(n+1)\sum_{k=0}^n\bi{n+1}{k+1}\bi{n+1+k}k\l(\f{b-1}2\r)^k-n\sum_{k=0}^{n-1}\bi n{k+1}\bi{n+k}k\l(\f{b-1}2\r)^k
\\=&\sum_{k=0}^n\((n+1+k)\bi{n+1}{k+1}-n\bi n{k+1}\)\bi{n+k}k\l(\f{b-1}2\r)^k
\\=&(2n+1)\sum_{k=0}^n\bi nk\bi{n+k}k\l(\f{b-1}2\r)^k=(2n+1)D_n(m)=(2n+1)T_n(b,c).
\endalign$$
Therefore, by the induction hypothesis, we have
$$\align&(n+1)\sum_{k=0}^n\bi{n+1}{k+1}\bi{n+1+k}k\l(\f{b-1}2\r)^k
\\=&\sum_{k=0}^{n-1}(2k+1)T_k(b,c)+(2n+1)T_n(b,c)=\sum_{k=0}^n(2k+1)T_k(b,c).
\endalign$$
This proves (1.15) with $n$ replaced by $n+1$.

 Let $p$ be a prime not dividing $b-1=2m$. It is easy to see that
$$\bi{2p-1}{p-1}=\prod_{j=1}^{p-1}\l(1+\f pj\r)\eq1+p\sum_{j=1}^{p-1}\f1j=1+p\sum_{j=1}^{(p-1)/2}\l(\f1j+\f1{p-j}\r)\eq1\ \ (\mo\ p^2).$$
In light of (1.15),
$$\align&\f1p\sum_{k=0}^{p-1}(2k+1)T_k(b,c)=\sum_{k=0}^{p-1}\bi p{k+1}\bi{p+k}km^k
\\=&\bi{2p-1}{p-1}m^{p-1}+\sum_{k=0}^{p-2}\bi p{k+1}\bi{p+k}km^k
\\\eq&m^{p-1}+\sum_{k=0}^{p-2}\bi p{k+1}m^k=m^{p-1}+\f{(m+1)^p-m^p-1}m
\\\eq&1+\f {(m+1)^p-(m+1)}m=1+\f{b+1}{b-1}\(\l(\f{b+1}2\r)^{p-1}-1\)\ (\mo\ p^2)
\endalign$$
and hence (1.16) follows.

Now we show (1.17). In view of Lemma 3.3,
$$\align\sum_{n=0}^{p-1}(2n+1)^2T_n(b,c)
=&\sum_{n=0}^{p-1}(2n+1)^2\sum_{k=0}^n\bi{n+k}{2k}\bi{2k}km^k
\\=&\sum_{k=0}^{p-1}\bi{2k}km^k\sum_{n=0}^{p-1}(2n+1)^2\bi{n+k}{2k}
\\=&(4p^2-1)\sum_{k=0}^{p-1}\f{p-k}{2k+3}\bi{p+k}{2k}\bi{2k}km^k
\\=&(4p^2-1)\sum_{k=0}^{p-1}\f{pm^k}{2k+3}\prod_{0<j\ls k}\l(\f{p^2}{j^2}-1\r)
\\\eq&-\sum_{k=0}^{p-1}\f{p(-m)^k}{2k+3}\pmod{p^2}
\\\eq&-(-m)^{(p-3)/2}\eq\f1m\l(\f{-m}p\r)\pmod{p}.
\endalign$$
This proves (1.17).
\qed

\heading{4. Proofs of Theorems 1.8 and 1.9}\endheading

\medskip

\noindent{\it Proof of Theorem} 1.8(i). We first prove (1.19)
by induction.

When $n=1$, both sides of (1.19) are equal to $b$.

 Now assume that (1.19) holds for a fixed integer $n\gs1$. Then
$$\align&b\sum_{k=0}^{(n+1)-1}(2k+1)T_k(b,c)^2(-d)^{(n+1)-1-k}
\\=&b(2n+1)T_n(b,c)^2-bd\sum_{k=0}^{n-1}(2k+1)T_k(b,c)^2(-d)^{n-1-k}
\\=&b(2n+1)T_n(b,c)^2-dnT_n(b,c)T_{n-1}(b,c)
\\=&(n+1)T_n(b,c)T_{n+1}(b,c).
\endalign$$
This concludes the induction step.

 Now we fix a positive integer $n$ and want to show (1.18).
As in the proof of Theorem 1.2(i),
$$T_n(b,c)\eq\cases \bi n{n/2}c^{n/2}\ (\mo\ b)&\t{if}\ 2\mid n,\\0\ (\mo\ b)&\t{if}\ 2\nmid n.\endcases$$

When $b\not=0$,  $b$ divides $T_n(b,c)$ or $T_{n-1}(b,c)$ since $n$
or $n-1$ is odd, therefore (1.18) follows from (1.19).

Now it remains to consider the case $b=0$. Note that $T_k(0,c)=0$
for $k=1,3,5,\ldots$, and $T_k(0,c)=\bi k{k/2}c^{k/2}$ for
$k=0,2,4,\ldots$. Thus
$$\align&\sum_{k=0}^{n-1}(2k+1)T_k(0,c)^2(4c-0^2)^{n-1-k}
\\=&\sum_{k=0}^{\lfloor(n-1)/2\rfloor}(4k+1)\l(\bi{2k}kc^k\r)^2(4c)^{n-1-2k}
\\=&(4c)^{n-1}\sum_{k=0}^{\lfloor(n-1)/2\rfloor}(4k+1)\f{\bi{2k}k^2}{16^k}.
\endalign$$
By induction, for any $m\in\N$ we have the identity
$$\sum_{k=0}^m(4k+1)\f{\bi{2k}k^2}{16^k}=\f{(m+1)^2}{16^m}\bi{2m+1}m^2=\f{(2m+1)^2}{16^m}\bi{2m}m^2,$$
which was pointed out to the author by R. Tauraso. It follows that
$$4^{n-1}\sum_{k=0}^{\lfloor(n-1)/2\rfloor}(4k+1)\f{\bi{2k}k^2}{16^k}=n^2\bi{n-1}{\lfloor n/2\rfloor}^2.$$
Therefore
$$\sum_{k=0}^{n-1}(2k+1)T_k(0,c)^2(4c-0^2)^{n-1-k}\eq0\ (\mo\ n^2)$$
and hence (1.18) holds when $b=0$. We are done. \qed

\proclaim{Lemma 4.1} Let $b,c\in\Z$ and $d=b^2-4c$. For any $n\in\N$ we have
$$T_n(b,c)^2=\sum_{k=0}^n\bi{n+k}{2k}\bi{2k}k^2c^kd^{n-k}.\tag4.1$$
\endproclaim
\Proof. If $d=0$ (i.e., $b^2=4c$), then
$$T_n(b,c)=[x^n]\l(x^2+bx+\f{b^2}4\r)^n=[x^n]\l(x+\f b2\r)^{2n}=\bi{2n}n\f{b^n}{2^n}$$
and hence (4.1) holds.

Now assume that $d\not=0$.  It is known that
$$\sum_{k=0}^n\bi{n+k}{2k}\bi{2k}k^2x^k(x+1)^k=\(\sum_{k=0}^n\bi{n}k\bi{n+k}kx^k\)^2$$
(cf. [S2, Lemma 3.2]), which is actually a special case of the famous Clausen identity for hypergeometric series.
Therefore
$$\align T_n(b,c)^2=&\((\sqrt d)^n P_n\l(\f b{\sqrt d}\r)\)^2=d^n D_n\l(\f{b/\sqrt d-1}2\r)^2
\\=&d^n\sum_{k=0}^n\bi{n+k}{2k}\bi{2k}k^2\l(\f{b/\sqrt d-1}2\r)^k\l(\f{b/\sqrt d+1}2\r)^k
\\=&d^n\sum_{k=0}^n\bi{n+k}{2k}\bi{2k}k^2\l(\f{b^2/d-1}4\r)^k=\sum_{k=0}^n\bi{n+k}{2k}\bi{2k}k^2c^kd^{n-k}.
\endalign$$
This completes the proof. \qed

\proclaim{Lemma 4.2} For any $k\in\N$ and $n\in\Z^+$ we have
$$\sum_{m=0}^{n-1}(2m+1)\bi{m+k}{2k}=\f{n(n-k)}{k+1}\bi{n+k}{2k}.\tag4.2$$
\endproclaim
\Proof. (4.2) can be easily proved by induction on $n$. \qed

\medskip
\noindent{\it Proof of Theorem} 1.8(ii). Let $n\in\Z^+$. In view of Lemmas 4.1 and 4.2, we have
$$\align&\sum_{m=0}^{n-1}(2m+1)T_m(b,c)^2d^{n-1-m}
\\=&\sum_{m=0}^{n-1}(2m+1)d^{n-1-m}\sum_{k=0}^m\bi{m+k}{2k}\bi{2k}k^2c^kd^{m-k}
\\=&\sum_{k=0}^{n-1}\bi{2k}k^2c^kd^{n-1-k}\sum_{m=0}^{n-1}(2m+1)\bi{m+k}{2k}
\\=&\sum_{k=0}^{n-1}\bi{2k}k^2c^kd^{n-1-k}\f{n(n-k)}{k+1}\bi{n+k}{2k}
\\=&n\sum_{k=0}^{n-1}(n-k)\bi nk\bi{n+k}kC_kc^kd^{n-1-k}
\\=&n^2\sum_{k=0}^{n-1}\bi{n-1}k\bi{n+k}kC_kc^kd^{n-1-k}.
\endalign$$
This proves (1.20).

Now assume $c\not=0$ and let $p$ be an odd prime not dividing $d$. By (1.20),
$$\f1{p^2}\sum_{k=0}^{p-1}(2k+1)\f{T_k(b,c)^2}{d^k}=\sum_{k=0}^{p-1}\bi{p-1}k\bi{p+k}kC_k\f{c^k}{d^k}.$$
For $k=0,1,\ldots,p-1$, clearly
$$\align\bi{p-1}k\bi{p+k}k=&\prod_{0<j\ls k}\l(\f{p-j}j\cdot\f{p+j}j\r)
=(-1)^k\prod_{0<j\ls k}\l(1-\f{p^2}{j^2}\r)
\\\eq&(-1)^k\l(1-p^2H_k^{(2)}\r)\pmod{p^4},
\endalign$$
where $H_k^{(2)}=\sum_{0<j\ls k}1/j^2$. Thus
$$\align\f1{p^2}\sum_{k=0}^{p-1}(2k+1)\f{T_k(b,c)^2}{d^k}
\eq&\sum_{k=0}^{p-1}C_k\l(-\f cd\r)^k(1-p^2H_k^{(2)})\pmod{p^4}
\\\eq&\sum_{k=0}^{p-1}C_k\l(-\f cd\r)^k\pmod{p^2}.
\endalign$$
If $p\mid c$, then $(\f dp)=(\f{b^2}p)=1$ and hence (1.21) follows.
In the case $p\nmid c$, we take an integer $m\eq-d/c\pmod{p^2}$ and then get
$$\f1{p^2}\sum_{k=0}^{p-1}(2k+1)\f{T_k(b,c)^2}{d^k}\eq\sum_{k=0}^{p-1}\f{C_k}{m^k}\pmod{p^2}.$$
By [Su3, Lemma 2.1],
$$\align\sum_{k=1}^{p-1}\f{C_k}{m^k}\eq&\f{m-4}2\l(1-\l(\f{m(m-4)}p\r)\r)
\\\eq&-\f{d+4c}{2c}\l(1-\l(\f{d(d+4c)}p\r)\r)=\f{b^2}{2c}\l(\l(\f dp\r)-1\r)\pmod p.
\endalign$$
(Moreover, the author [Su1] determined $\sum_{k=1}^{p-1}C_k/m^k$ mod $p^2$ in terms of Lucas sequences.)
So (1.21) is valid. We are done. \qed

\Remark\ 4.3. Let $p>3$ be a prime. As $D_k=T_k(3,2)$, by refining the proof of Theorem 1.8(ii) and using two auxiliary congruences
$$\sum_{k=1}^{p-1}(-2)^kC_k\eq-4p\,q_p(2)\pmod{p^3}$$
and $$\sum_{k=1}^{p-1}(-2)^kC_kH_k^{(2)}\eq2q_p(2)^2\pmod p$$
(the author has a proof of them), we get
$$\sum_{k=0}^{p-1}(2k+1)D_k^2\eq p^2-4p^3q_p(2)-2p^4q_p(2)^2\ (\mo\ p^5).$$

\proclaim{Lemma 4.4} Let $b,c\in\Z$. Suppose that $p>3$ is a prime not dividing $d=b^2-4c$.
Then
$$\sum_{k=0}^{p-1}\f{T_k(b,c)^2}{d^k}\eq\l(\f{16c}d\r)^{(p-1)/2}+p\sum\Sb k=0\\k\not=(p-1)/2\endSb^{p-1}\f{\bi{2k}k}{2k+1}\l(-\f cd\r)^k
\pmod{p^3}.\tag4.3$$
\endproclaim
\Proof. With the help of (4.1), we have
$$\align\sum_{n=0}^{p-1}\f{T_n(b,c)^2}{d^n}=&\sum_{n=0}^{p-1}\f1{d^n}\sum_{k=0}^n\bi{n+k}{2k}\bi{2k}k^2c^kd^{n-k}
\\=&\sum_{k=0}^{p-1}\bi{2k}k^2\f{c^k}{d^k}\sum_{n=k}^{p-1}\bi{n+k}{2k}=\sum_{k=0}^{p-1}\bi{2k}k^2\bi{p+k}{2k+1}\l(\f cd\r)^k
\\=&\sum_{k=0}^{p-1}\f p{2k+1}\bi{2k}k\(\prod_{0<j\ls k}\f{p^2-j^2}{j^2}\)\l(\f cd\r)^k
\endalign$$
and hence
$$\align\sum_{n=0}^{p-1}\f{T_n(b,c)^2}{d^n}
\eq&\sum_{k=0}^{p-1}\f{p(-1)^k}{2k+1}\bi{2k}k\l(1-p^2H_k^{(2)}\r)\l(\f cd\r)^k\pmod{p^4}
\\\eq&(-1)^{(p-1)/2}\bi{p-1}{(p-1)/2}\l(1-p^2H_{(p-1)/2}^{(2)}\r)\l(\f cd\r)^{(p-1)/2}
\\&+p\sum\Sb k=0\\k\not=(p-1)/2\endSb^{p-1}\f{\bi{2k}k}{2k+1}\l(-\f cd\r)^k\pmod{p^3}.
\endalign$$
As Wolstenholme observed, $H_{p-1}^{(2)}\eq0\pmod p$ since $\sum_{j=1}^{p-1}1/(2j)^2\eq\sum_{k=1}^{p-1}1/k^2\pmod p$.
Therefore
$$H_{(p-1)/2}^{(2)}\eq\f12\sum_{k=1}^{(p-1)/2}\l(\f1{k^2}+\f1{(p-k)^2}\r)=\f{H_{p-1}^{(2)}}2\eq0\pmod p.$$
Recall Morley's congruence (cf. [M])
$$\bi{p-1}{(p-1)/2}\eq(-1)^{(p-1)/2}4^{p-1}\pmod{p^3}.$$
So we have
$$(-1)^{(p-1)/2}\bi{p-1}{(p-1)/2}\l(1-p^2H_{(p-1)/2}^{(2)}\r)\eq4^{p-1}\pmod{p^3}$$
and hence (4.3) follows. \qed

\medskip
\noindent{\it Proof of Theorem} 1.9. (i) Applying Lemma 4.4 with $b=6$ and $c=-3$ we get
$$\sum_{k=0}^{p-1}\f{T_k(6,-3)^2}{48^k}\eq\l(\f{-1}p\r)+p\sum^{p-1}\Sb k=0\\k\not=(p-1)/2\endSb\f{\bi{2k}k}{(2k+1)16^k}\pmod {p^3}.$$
By [Su2, (1.4)-(1.5)],
$$\sum_{k=0}^{(p-3)/2}\f{\bi{2k}k}{(2k+1)16^k}\eq0\pmod{p^2}$$ and
$$\sum_{k=(p+1)/2}^{p-1}\f{\bi{2k}k}{(2k+1)16^k}\eq\f p3E_{p-3}\pmod{p^2}.$$
So (1.22) follows.
\medskip

(ii) Now we prove (1.23) and (1.24). Since $p\mid\bi{2k}k$ for every $k=(p+1)/2,\ldots,p-1$, by Lemma 4.4 with $b=2$ and $c\in\{-1,-3\}$ we obtain
$$\sum_{k=0}^{p-1}\f{T_k(2,-1)^2}{8^k}\eq(-2)^{(p-1)/2}+p\sum_{k=0}^{(p-3)/2}\f{\bi{2k}k}{(2k+1)8^k}\pmod{p^2}\tag4.4$$
and
$$\sum_{k=0}^{p-1}\f{T_k(2,-3)^2}{16^k}\eq(-3)^{(p-1)/2}+p\sum_{k=0}^{(p-3)/2}\f{\bi{2k}k}{(2k+1)}\l(\f 3{16}\r)^k\ (\mo\ p^2).\tag4.5$$

For $n\in\N$ define
$$u_n=(2n+1)\sum_{k=0}^n\bi{n+k}{2k}\f{(-2)^k}{2k+1}
\ \ \t{and}\ \ v_n=(2n+1)\sum_{k=0}^n\bi{n+k}{2k}\f{(-3)^k}{2k+1}.$$
Via the Zeilberger algorithm (cf. [PWZ]) we find the recurrence relations
$$u_n+u_{n+2}=0\ \quad\t{and}\quad\ v_n+v_{n+1}+v_{n+2}=0.$$
So, by induction we have
$$u_n=(-1)^{n(n-1)/2}=\l(\f{-2}{2n+1}\r)\ \ \t{and}\ \ v_n=\l(\f{2n+1}3\r)$$
for every $n=0,1,2,\ldots$.
Taking $n=(p-1)/2$ and noting that
$\bi{n+k}{2k}\eq\bi{2k}k/(-16)^k\pmod{p^2}$ for $k=0,\ldots,n$ (cf. [S1, Lemma 2.2]), we then obtain
$$(-2)^{(p-1)/2}+p\sum_{k=0}^{(p-3)/2}\f{\bi{2k}k}{(2k+1)8^k}\eq u_{(p-1)/2}=\l(\f{-2}p\r)\pmod {p^3}$$
and
$$(-3)^{(p-1)/2}+p\sum_{k=0}^{(p-3)/2}\f{\bi{2k}k}{(2k+1)}\l(\f 3{16}\r)^k\eq v_{(p-1)/2}=\l(\f p3\r)\pmod {p^3}.$$
Combining these with (4.4) and (4.5) we immediately get (1.23) and (1.24).
\medskip

(iii) Finally we show (1.25). Applying (4.1) with $b=3$ and $c=2$ we obtain
$$\sum_{k=0}^n\bi{n+k}{2k}\bi{2k}k^22^k=D_n^2.$$
Therefore
$$\align\sum_{n=1}^{p-1}\f{D_n^2-1}{n^2}=&\sum_{n=1}^{p-1}\f1{n^2}\sum_{k=1}^n\bi{n+k}{2k}\bi{2k}k^22^k
\\=&\sum_{k=1}^{p-1}2^k\bi{2k}k^2\ \sum_{n=k}^{p-1}\f{\bi{n+k}{2k}}{n^2}
=\sum_{k=1}^{p-1}2^k\bi{2k}k^2\sum_{r=0}^{p-1-k}\f{\bi{2k+r}{r}}{(k+r)^2}.
\endalign$$
If $k\in\{(p+1)/2,\ldots,p-1\}$ then $p\mid\bi{2k}k$.
For each $k=1,\ldots,(p-1)/2$, clearly
$$\sum_{r=0}^{p-1-k}\f{\bi{2k+r}r}{(k+r)^2}=4\sum_{r=0}^{p-1-k}\f{(-1)^r\bi{-2k-1}r}{(-2k-2r)^2}
\eq4\sum_{r=0}^{p-1-2k}\f{(-1)^r\bi{p-1-2k}r}{(p-2k-2r)^2}\pmod p.$$
By [Su2, (2.5)], we have the identity
$$\sum_{r=0}^{2n}\f{(-1)^r\bi{2n}r}{(2n+1-2r)^2}=\f{(-16)^n}{(2n+1)^2\bi{2n}n}.$$
Also, $H_{p-1}^{(2)}=\sum_{k=1}^{p-1}1/k^2\eq0\pmod p$.
So, by the above, we have
$$\align\sum_{k=1}^{p-1}\f{D_k^2}{k^2}\eq&\sum_{k=1}^{(p-1)/2}2^k\bi{2k}k^2\f{4(-16)^{(p-1)/2-k}}
{(p-2k)^2\bi{p-1-2k}{(p-1)/2-k}}
\\\eq&\sum_{k=1}^{(p-1)/2}\f{2^k\bi{2k}k^24^{(p-1)/2-k}}{k^2\bi{(p-1-2k}{(p-1)/2-k}/(-4)^{(p-1)/2-k}}\pmod{p}.\endalign$$
For each $k\in\{1,\ldots,(p-1)/2\}$, obviously
$$\align\f{\bi{2k}k}{(-4)^k}=\bi{-1/2}k\eq&\bi{(p-1)/2}k=\bi{(p-1)/2}{(p-1)/2-k}
\\\eq&\bi{-1/2}{(p-1)/2-k}=\f{\bi{p-1-2k}{(p-1)/2-k}}{(-4)^{(p-1)/2-k}}\pmod p.
\endalign$$
Therefore
$$\sum_{k=1}^{p-1}\f{D_k^2}{k^2}\eq\sum_{k=1}^{(p-1)/2}\f{2^k\bi{2k}k^22^{p-1}/4^k}{k^2\bi{2k}k/(-4)^k}
\eq\sum_{k=1}^{(p-1)/2}\f{(-2)^k\bi{2k}k}{k^2}\pmod{p}$$
and hence
$$\sum_{k=1}^{p-1}\f{D_k^2}{k^2}\eq\sum_{k=1}^{p-1}\f{(-2)^k}{k^2}\bi{2k}k\ \pmod p.\tag4.6$$

Let
$$v_n=2^n+2^{-n}\quad\t{and}\quad w_n=(-1)^n+2^{-n}\quad\t{for all}\ n\in\N.$$
It is easy to see that
$$v_{n+1}=\f 52v_n-v_{n-1}\quad\t{and}\quad w_{n+1}=-\f12w_n+\f12w_{n-1}\quad\t{for all}\ n\in\Z^+.$$
Thus, applying [MT, (42)] with $t=-1/2$ we obtain
$$\align-\f14\sum_{k=1}^{p-1}\f{(-2)^k}{k^2}\bi{2k}k&\eq\f{v_p+2w_p+2^{-p}-2}{p^2}+\sum_{k=1}^{p-1}\f{v_k}{k^2}
\\&=2^{-p}\l(\f{2^p-2}p\r)^2+\sum_{k=1}^{p-1}\f{2^k}{k^2}+\sum_{k=1}^{p-1}\f{2^{-(p-k)}}{(p-k)^2}
\\&\eq2q_p(2)^2+\f32\sum_{k=1}^{p-1}\f{2^k}{k^2}\ \pmod{p}.
\endalign$$
Recall that $\sum_{k=1}^{p-1}2^k/k^2\eq-q_p(2)^2\ (\mo\ p)$ (which was conjectured by L. Skula and proved by A. Granville [Gr]). So we have
$$\sum_{k=1}^{p-1}\f{(-2)^k}{k^2}\bi{2k}k\eq-2q_p(2)^2\ \pmod p.\tag4.7$$

Combining (4.6) and (4.7) we finally get (1.25).
This ends the proof. \qed

\heading{5. More conjectures for further research}\endheading

Motivated by part (ii) of Theorem 1.5, we raise the following conjecture.

\proclaim{Conjecture 5.1} Let $x$ be any integer.  Then
$$\sum_{k=0}^{n-1}(2k+1)D_k(x)^m\eq0\ \ (\mo\ n)$$
for all $m,n\in\Z^+$. If $p$ is a prime not dividing $x(x+1)$, then
$$\sum_{k=0}^{p-1}(2k+1)D_k(x)^3\eq p\l(\f{-4x-3}p\r)\ (\mo\ p^2)$$
and
$$\sum_{k=0}^{p-1}(2k+1)D_k(x)^4\eq p\ (\mo\ p^2).$$
\endproclaim

Now we propose the following conjecture related to Theorem 1.2(ii).

\proclaim{Conjecture 5.2} Let $b,c\in\Z$. For any $n\in\Z^+$ we have
$$\sum_{k=0}^{n-1}(8ck+4c+b)T_k(b,c^2)^2(b-2c)^{2(n-1-k)}\eq0\ (\mo\ n).$$
If $p$ is an odd prime not dividing $b(b-2c)$, then
$$\sum_{k=0}^{p-1}(8ck+4c+b)\f{T_k(b,c^2)^2}{(b-2c)^{2k}}\eq p(b+2c)\l(\f{b^2-4c^2}p\r)\ (\mo\ p^2).$$
\endproclaim
\Remark\ 5.3. Conjecture 5.2 in the case $b=c=1$ yields the first part of Conjecture 1.1.
\medskip

By Theorem 1.2(ii), if $p$ is an odd prime then
$$\sum_{k=0}^{p-1}\f{T_k(4,1)^2}{2^{2k}}\eq\sum_{k=0}^{p-1}\f{T_k(4,1)^2}{6^{2k}}\eq\l(\f{-1}p\r)\ (\mo\ p).$$
Motivated by this and (1.22)-(1.24), we now give a further conjecture.
\proclaim{Conjecture 5.4} Let $p$ be an odd prime.
We have
$$\sum_{k=0}^{p-1}\f{T_k(2,2)^2}{4^k}\eq\sum_{k=0}^{p-1}\f{\bi{2k}k^2}{8^k}
\ \ \l(\mo\ p^{(5+(\f{-1}p))/2}\r).$$
If $p>3$, then
$$\sum_{k=0}^{p-1}\f{T_k(4,1)^2}{4^k}\eq\sum_{k=0}^{p-1}\f{T_k(4,1)^2}{36^k}\eq\l(\f{-1}p\r)\ (\mo\ p^2).$$
\endproclaim

Now we raise a conjecture related to Theorem 1.2(iii).
\proclaim{Conjecture 5.5} Let $b,c\in\Z$ and $d=b^2-4c$. For any $n\in\Z^+$ we have
$$\sum_{k=0}^{n-1}T_k(b,c)M_k(b,c)d^{n-1-k}\eq0\ (\mo\ n).$$
If $p$ is an odd prime not dividing $cd$, then
$$\sum_{k=0}^{p-1}\f{T_k(b,c)M_k(b,c)}{d^k}\eq\f{pb^2}{2c}\(\l(\f{d}p\r)-1\)\ (\mo\ p^2).$$
\endproclaim

By Conjecture 5.5, for any prime $p>3$ we should have
$$\sum_{k=0}^{p-1}\f{T_k(3,3)M_k(3,3)}{(-3)^k}\eq\f{3p}2\l(\l(\f p3\r)-1\r)\ (\mo\ p^2).$$
This can be further strengthened.

\proclaim{Conjecture 5.6} Let $p>3$ be a prime. Then
$$\sum_{k=0}^{p-1}\f{T_k(3,3)M_k(3,3)}{(-3)^k}\eq\cases2p^2\ (\mo\ p^3)&\t{if}\ p\eq1\ (\mo\ 3),
\\p^3-p^2-3p\ (\mo\ p^4)&\t{if}\ p\eq2\ (\mo\ 3).\endcases$$
\endproclaim

In view of Theorem 1.2(ii), for $b,c\in\Z$ and a prime $p\nmid (b-2c)$, it is natural to
investigate whether the sum $\sum_{k=0}^{p-1}T_k(b,c^2)^3/(b-2c)^{3k}$ mod $p$
has a pattern. This leads us to raise the following two conjectures.

\proclaim{Conjecture 5.7} Let $p>3$ be a prime. Then
$$\align&\l(\f 3p\r)\sum_{k=0}^{p-1}\f{T_k(2,3)^3}{8^k}\eq\sum_{k=0}^{p-1}\f{T_k(2,3)^3}{(-64)^k}
\\\eq&\sum_{k=0}^{p-1}\f{T_k(2,9)^3}{(-64)^k}\eq\l(\f 3p\r)\sum_{k=0}^{p-1}\f{T_k(2,9)^3}{512^k}
\\\eq&\cases4x^2-2p\ (\mo\ p^2)&\t{if}\ p\eq1,7\ (\mo\ 24)\ \t{and}\ p=x^2+6y^2,
\\2p-8x^2\ (\mo\ p^2)&\t{if}\ p\eq5,11\ (\mo\ 24)\ \t{and}\ p=2x^2+3y^2,
\\0\ (\mo\ p^2)&\t{if}\ (\f {-6}p)=-1.\endcases
\endalign$$
And
$$\align \sum_{k=0}^{p-1}(3k+2)\f{T_k(2,3)^3}{8^k}\eq&p\l(3\l(\f 3p\r)-1\r)\pmod{p^2},
\\\sum_{k=0}^{p-1}(3k+1)\f{T_k(2,3)^3}{(-64)^k}\eq&p\l(\f{-2}p\r)\pmod{p^3}.
\endalign$$
When $(\f{-6}p)=1$ we have
$$\sum_{k=0}^{p-1}(72k+47)\f{T_k(2,9)^3}{(-64)^k}\eq 42p\ \ (\mo\ p^2)$$
and
$$\sum_{k=0}^{p-1}(72k+25)\f{T_k(2,9)^3}{512^k}\eq 12p\l(\f 3p\r)\ (\mo\ p^2).$$
Also,
$$\sum_{k=0}^{n-1}(3k+2)T_k(2,3)^38^{n-1-k}\eq0\ \ (\mo\ 2n)$$
and
$$\sum_{k=0}^{n-1}(3k+1)T_k(2,3)^3(-64)^{n-1-k}\eq0\ \ (\mo\ n)$$
for every positive integer $n$.
\endproclaim
\Remark\ 5.8. Let $p>3$ be a prime. If $p\eq1,7\ (\mo\
24)$ then $p=x^2+6y^2$ for some $x,y\in\Z$; if $p\eq5,11\ (\mo\ 24)$ then
$p=2x^2+3y^2$ for some $x,y\in\Z$. The reader may consult [BEW] and [Co] for such known facts.

\medskip

\proclaim{Conjecture 5.9} Let $p>3$ be a prime. Then
$$\align&\l(\f {2}p\r)\sum_{k=0}^{p-1}\f{T_k(18,49)^3}{8^{3k}}\eq\sum_{k=0}^{p-1}\f{T_k(18,49)^3}{16^{3k}}
\\\eq&\cases 4x^2-2p\pmod{p^2}&\t{if}\ p\eq1\pmod4\ \&\ p=x^2+y^2\ (2\nmid x,\ 2\mid y),
\\0\pmod{p^2}&\t{if}\ p\eq3\pmod4.\endcases
\endalign$$
And
$$\align&\l(\f {-1}p\r)\sum_{k=0}^{p-1}\f{T_k(10,49)^3}{(-8)^{3k}}\eq\l(\f 6p\r)\sum_{k=0}^{p-1}\f{T_k(10,49)^3}{12^{3k}}
\\\eq&\cases 4x^2-2p\pmod{p^2}&\t{if}\ p\eq1,3\pmod8\ \&\ p=x^2+2y^2\ (x,y\in\Z),
\\0\pmod{p^2}&\t{if}\ (\f{-2}p)=-1,\ \t{i.e.},\ p\eq5,7\pmod8.\endcases
\endalign$$
Also,
$$\align\sum_{k=0}^{p-1}(7k+4)\f{T_k(10,49)^3}{(-8)^{3k}}\eq&\f p{14}\l(\f 2p\r)\l(65-9\l(\f p3\r)\r)\pmod{p^2},
\\\sum_{k=0}^{p-1}(7k+3)\f{T_k(10,49)^3}{12^{3k}}\eq&\f{3p}{28}\l(13+15\l(\f p3\r)\r)\pmod{p^2}.
\endalign$$
For each $n=1,2,3,\ldots$ we have
$$\sum_{k=0}^{n-1}(7k+4)T_k(10,49)^3(-8^3)^{n-1-k}\eq0\pmod{4n}$$
and
$$\sum_{k=0}^{n-1}(7k+3)T_k(10,49)^3(12^3)^{n-1-k}\eq0\pmod{n}.$$
\endproclaim

\enddocument